\documentclass[12pt]{article}
\usepackage{amssymb}
\usepackage{amsmath}
\allowdisplaybreaks

\newcommand{\di}{\partial}

\numberwithin{equation}{section}

\newcommand{\rom}[1]{{\rm #1}}

\newtheorem{theorem}{Theorem}[section]

\newtheorem{lemma}{Lemma}[section]
\newtheorem{remark}{Remark}[section]
\newtheorem{corollary}{Corollary}[section]
\newtheorem{proposition}{Proposition}[section]

\newcommand{\R}{\mathbb{R}}
\newcommand{\Z}{\mathbb{Z}}
\newcommand{\N}{\mathbb{N}}
\newcommand{\eps}{\varepsilon}
\newcommand{\la}{\langle}
\newcommand{\ra}{\rangle}

\begin{document}

\begin{center}{\Large \bf
 Diffusion approximation for equilibrium Kawasaki dynamics in
 continuum}\end{center}

{\large Yuri G. Kondratiev}\\
 Fakult\"at f\"ur Mathematik, Universit\"at
Bielefeld, Postfach 10 01 31, D-33501 Bielefeld, Germany; BiBoS,
Univ.\ Bielefeld, Germany; Kiev-Mohyla Academy, Kiev, Ukraine\\
 e-mail:
\texttt{kondrat@mathematik.uni-bielefeld.de}\vspace{2mm}

{\large Oleksandr V.  Kutoviy}\\
 Fakult\"at f\"ur Mathematik, Universit\"at Bielefeld, Postfach 10
01 31, D-33501 Bielefeld, Germany; BiBoS,
Univ.\ Bielefeld, Germany\\
e.mail: \texttt{kutoviy@mathematik.uni-bielefeld.de}\vspace{2mm}

{\large Eugene W. Lytvynov}\\
Department of Mathematics, University of Wales Swansea, Singleton
Park, Swansea SA2 8PP, U.K.\\ e-mail:
\texttt{e.lytvynov@swansea.ac.uk}

{\small
\begin{center}
{\bf Abstract}
\end{center}
\noindent A Kawasaki dynamics in continuum is a dynamics of an
infinite system of interacting particles in $\R^d$ which randomly
hop over the space. In this paper, we deal with an equilibrium
Kawasaki dynamics which has a Gibbs measure $\mu$ as invariant
measure. We study a diffusive  limit of such a dynamics, derived
through a scaling of both the jump rate and time. Under weak
assumptions on the potential of pair interaction, $\phi$, (in
particular, admitting a singularity of $\phi$ at zero), we prove
that, on a set of smooth local functions,  the generator of the
scaled dynamics converges  to the generator of the gradient stochastic dynamics.
If the set on which the generators converge is a core for the
diffusion generator, the latter result implies the weak convergence
of finite-dimensional distributions of the corresponding equilibrium
processes. In particular, if the potential $\phi$ is from
$C_{\mathrm b}^3(\R^d)$ and sufficiently quickly converges to zero
at infinity, we conclude the
convergence of the processes from a result in [Choi~{\it et~al.}, {
J. Math.\ Phys.}\ 39 (1998) 6509--6536].
 }
\vspace{3mm}

\vspace{1.5mm}

\vspace{1.5mm}

\section{Introduction} A Kawasaki dynamics in continuum is a dynamics of an infinite system
of interacting particles in $\R^d$ which randomly hop over the
space. The generator of such a dynamics has the form
\begin{equation}\label{hjfyt} (HF)(\gamma)=-\sum_{x\in\gamma}\int_{\R^d}dy\, c
(\gamma,x,y)
(F(\gamma\setminus x\cup y)-F(\gamma)),\quad
\gamma\in\Gamma. \end{equation}
Here,
 $\Gamma$ denotes the
configuration space over $\R^d$, i.e., the space of all locally
finite subsets of $\R^d$, and, for simplicity of notations, we just
write $x$ instead of $\{x\}$. The coefficient $c(\gamma,x,y)$
describes the rate at which the particle $x$ of the configuration
$\gamma$ jumps to $y$.

Let $\mu$ denote a Gibbs measure on $\Gamma$ which
corresponds to an activity parameter $z>0$ and a
potential of pair interaction $\phi$.
In this paper, we will deal with an   equilibrium Kawasaki dynamics
which has $\mu$ as invariant measure. More precisely, 
we will consider an equilibrium Kawasaki dynamics whose   generator \eqref{hjfyt} has the coefficient
$c(\gamma,x,y)$ of the form
\begin{equation}\label{iuiy} c(\gamma,x,y)=a(x-y)
\exp\big[(1/2)E(x,\gamma\setminus x)-(1/2)E(y,\gamma\setminus x)\big]. \end{equation}
Here, for any   $\gamma\in\Gamma$ and  $u\in\R^d\setminus\gamma$, $E(u,\gamma)$ denotes the relative energy of interaction between the  particle at $u$ and  the configuration $\gamma$.
About the function $a(\cdot)$ in
\eqref{iuiy} we assume that it is non-negative, bounded, has a
compact support,  and $a(x)$ only depends on $|x|$.

Equation \eqref{iuiy} allows the following physical interpretation: particles from  $\gamma$ which have
a high relative energy of interaction with the rest of the configuration tend to jump to places where this relative energy will be low, i.e., particles tend to jump from high energy regions to low energy regions.

Note also that the bilinear (Dirichlet) form corresponding to the generator \eqref{hjfyt}, \eqref{iuiy} admits the following representation:
\begin{multline*} \mathcal E(F,G)=\frac z2\int_\Gamma\mu(d\gamma)
\int_{\R^d}dx\int_{\R^d}dy\, a(x-y)\exp\big[
-(1/2)E(x,\gamma)-(1/2)E(y,\gamma)\big]\\
\times (F(\gamma \cup y)-F(\gamma\cup x))
(G(\gamma\cup y)-G(\gamma\cup x)).\end{multline*}

Under very mild assumptions on the Gibbs measure $\mu$, it was
proved in \cite{KLR} that there indeed exists a Markov process on $\Gamma$
with {\it c\'adl\'ag} paths whose generator is given by \eqref{hjfyt}, \eqref{iuiy}. We assume that the initial distribution of this
dynamics is $\mu$, and perform a diffusive scaling of this dynamics.
More precisely, for each $\epsilon>0$, we consider the equilibrium
Kawasaki dynamics whose jump rate is given by formula \eqref{iuiy} in
which $a(\cdot)$ is replaced with the function
\begin{equation}\label{rtsra}
a_\epsilon(\cdot):=\epsilon^{-d}a(\cdot/\epsilon),\end{equation} and
we additionally scale time, multiplying it by $\epsilon^{-2}$. We
denote the generator of the obtained  dynamics by
$H^{(\epsilon)}$.

So, the aim of the paper is to show that the scaled dynamics
converges, as $\epsilon\to0$, to a diffusive dynamics on the
configuration space $\Gamma$. Our main result is that,
under weak assumptions on the pair potential  $\phi$ (in particular, we allow $\phi$
to have a singularity at zero),    the generator of the scaled dynamics,
$H^{(\epsilon)}$, converges, on a set of smooth
local functions,  to the generator of the (infinite-dimensional) {\it gradient stochastic
dynamics} (also called {\it interacting Brownian particles}), see e.g.\
\cite{AKR4,Fritz,GKLR,Lang77,Osa96,osada,osada1,Spohn,Yos96} and
the references therein. So, the limiting diffusive generator acts as follows:
\begin{align*}(
&H^{\mathrm{(dif)}}F)(\{x_k\}_{k=1}^\infty)\\
&\quad=\frac{c}2\sum_{i=1}^\infty \bigg(-
\Delta_{x_i}F(\{x_k\}_{k=1}^\infty)+\sum_{j\ne i}\la \nabla_{x_i}
F(\{x_k\}_{k=1}^\infty),\nabla\phi(x_i-x_j)\ra\bigg),\label{bhjgh}
\end{align*}
  where the constant $c$ is defined in  the equation \eqref{dtse}
  below. The corresponding stochastic process
informally solves
the following system of stochastic differential equations:
\begin{align*}
&dx_i(t)=-\frac{c}2\sum_{j\ne
i}\nabla\phi(x_i(t)-x_j(t))\,dt+\sqrt{c}\,dB_i(t),\quad
i\in\N,\\
& \{x_i(0)\}_{i=1}^\infty=\gamma\in\Gamma,
\end{align*}
where $(B_i)_{i=1}^\infty$ is a sequence of independent Brownian
motions.

If the set on which the generators converge is a core for the
diffusive generator $H^{\mathrm{(dif)}}$, then our main result implies
 the weak convergence of finite-dimensional
distributions of the corresponding equilibrium processes. In
particular,  if the potential $\phi$ is from $C_{\mathrm b}^3(\R^d)$
(hence, $\phi$ has no singularity at zero) and sufficiently
quickly converges to zero at infinity, then we conclude the convergence of the process from a
result by Choi~{\it et~al.} \cite{CPY}.

The paper is organized as follows. In Section \ref{tydey}, we recall  some basic facts of analysis on the configuration space $\Gamma$. In Section~\ref{wqwwwq}, we recall  conditions which guarantee the existence of a Gibbs measure on the configuration space.  In Sections~\ref{tyedtwe} and \ref{ytrtyed}, we
recall construction of the equilibrium Kawasaki dynamics in continuum, and the gradient stochastic dynamics, respectively.   In Section~\ref{yfyyty}, we formulate our main results. Finally, in Section~\ref{jhgftfru}, we present the proofs.

 \section{$K$-transform and correlation functions}
\label{tydey}

  The configuration space over $\R^d$, $d\in\N$,
is defined as the set of all subsets of $\R^d$ which are locally
finite: $$\Gamma:=\Gamma_{\R^d}:=\big\{\,\gamma\subset \R^d\mid
|\gamma_\Lambda|<\infty\text{ for each  }\Lambda\in \mathcal O_{\mathrm c}(\R^d)\,\big\}.$$
Here $|\cdot|$ denotes the cardinality of a set,
$\gamma_\Lambda:= \gamma\cap\Lambda$, and $\mathcal O_{\mathrm c}(\R^d)$ denotes the set of all open, relatively
compact subsets of $\R^d$. One can identify any
$\gamma\in\Gamma$ with the positive Radon measure
$\sum_{x\in\gamma}\eps_x\in{\cal M}(\R^d)$, where  $\eps_x$ is
the Dirac measure with mass at $x$, and  ${\cal M}(\R^d)$
 stands for the set of all
positive  Radon  measures on the Borel $\sigma$-algebra ${\cal
B}(\R^d)$. The space $\Gamma$ can be endowed with the relative
topology as a subset of the space ${\cal M}(\R^d)$ with the vague
topology, i.e., the weakest topology on $\Gamma$ with respect to
which  all maps $$\Gamma\ni\gamma\mapsto\la f,\gamma\ra:=\int_{\R^d}
f(x)\,\gamma(dx) =\sum_{x\in\gamma}f(x),\qquad f\in C_0(\R^d),$$ are
continuous. Here, $C_0(\R^d)$ is the space of all
continuous  functions on $\R^d$ with
compact support. We will denote by ${\cal B}(\Gamma)$ the Borel
$\sigma$-algebra on $\Gamma$.

Next, denote by $\Gamma_0$ the space of finite configurations in ${\mathbb
R}^d$:
\begin{eqnarray*}
\Gamma_0 := \bigsqcup_{n=0}^\infty \Gamma_0^{(n)},\quad
\Gamma^{(0)}_0 := \{\varnothing\},\quad \Gamma_0^{(n)}:=\{\eta
\subset {\mathbb R}^d \mid |\eta|= n \}, \quad n \in {\mathbb N}.
\end{eqnarray*}
Evidently, $\Gamma_0\subset\Gamma$.

Let $$\widetilde{({\mathbb R}^d)^n}= \big\{\,(x_1, \dots, x_n) \in
({\mathbb R}^{d })^n  \mid x_i \ne x_j \text{ for }  i \ne
j \, \big\}.$$ Let $S^n$ be the group of all permutations of
$\{1,\dots,n\}$ which  acts on $\widetilde{({\mathbb
R}^d)^n}$ by permuting the coordinates. Through the natural bijection
\begin{eqnarray}\label{eq366}
\widetilde{(\mathbb R^d)^n}/S^n \longleftrightarrow \Gamma_0^{(n)}
\end{eqnarray}
one defines a topology on $\Gamma^{(n)}_0$. The space $\Gamma_0$ is then
equipped  with the topology of disjoint union. Let ${\mathcal
B}(\Gamma_0)$ denote the Borel $\sigma$-algebra on $\Gamma_0$. It can be shown (see e.g.\ \cite{Kal})
that ${\mathcal B}(\Gamma_0)$ coincides with the trace $\sigma$-algebra of $\mathcal B(\Gamma)$ on $\Gamma_0$.
Note also
that each function $k:\Gamma_0\to\R$ may be identified   with the
sequence $(k^{(n)})_{n=0}^\infty$, where $k^{(0)}:=k(\{\varnothing\})$ and, for each $n\in\mathbb N$,
$k^{(n)}:\widetilde{(\mathbb R^d)^n}\to\R$ is a measurable,
symmetric function.


For any $\gamma \in {\Gamma}$, let $\sum_{\eta \Subset \gamma}$
denote the summation over all $\eta \subset \gamma$ such that
$\eta\in\Gamma_0$. For a function $G: \Gamma_0 \to {\mathbb R}$, the
$K$-transform of $G$ is defined by
\begin{eqnarray}\label{eq9}
(KG)(\gamma):= \sum_{\eta \Subset \gamma} G(\eta)
\end{eqnarray}
for each $\gamma \in \Gamma$ such that at least one of the series
$\sum_{\eta \Subset \gamma} G^+(\eta)$,  $\sum_{\eta \Subset
\gamma} G^-(\eta)$ converges. Here $G^{+} := \max \{ 0, G\}$ and
$G^{-} := \max \{ 0, -G\}$.

Let us fix a probability measure $\mu$ on $(\Gamma,{\mathcal
B}(\Gamma))$. The correlation measure of $\mu$ is
defined by
\begin{eqnarray*}
\rho_\mu(A) := \int_{\Gamma}(K\chi_A)(\gamma) \,\mu(d\gamma), \qquad
A \in {\mathcal B}(\Gamma_0),
\end{eqnarray*}
where $\chi_A$ denotes the indicator of the set $A$. The  $\rho_\mu$
is a measure on $(\Gamma_0, {\mathcal B} (\Gamma_0))$ (see \cite{KK}
for details, in particular, measurability issues).
Note that $\rho_\mu(\{\varnothing\})=1$.

The following proposition was proved in \cite{KK}, see also \cite{Len75a,Len75b}

\begin{proposition}\label{dftrert}
Let $G \in L^1(\Gamma_0, \rho_\mu)$\rom, then
$KG \in
L^1(\Gamma, \mu)$, the series in \eqref{eq9} is absolutely convergent for $\mu$-a.e.~$\gamma \in
\Gamma$, and
 $$ \| KG \|_{L^1(\mu)} \le
\| K|G|\, \|_{L^1(\mu)} = \| G \|_{L^1(\rho_\mu)}.$$  Moreover, then
\begin{eqnarray}\label{eq302}
\int_{\Gamma_0} G(\eta) \,\rho_\mu(d\eta) =
\int_{\Gamma}(KG)(\gamma) \,\mu(d\gamma).
\end{eqnarray}
\end{proposition}

The Lebesgue--Poisson measure $\lambda$ on $(\Gamma_0,{\mathcal
B}(\Gamma_0))$ is defined by
\begin{eqnarray*}
\lambda := \varepsilon_\varnothing+\sum_{n=1}^\infty
\frac{1}{n!}\,dx^{\otimes n},
\end{eqnarray*}
where $dx^{\otimes n}$ is defined via the bijection (\ref{eq366}).
Assume that the correlation measure $\rho_\mu$  is absolutely
continuous with respect to the Lebesgue--Poisson measure $\lambda$.
Denote $k_\mu:= d\rho_\mu/d\lambda$. Then the corresponding
functions $(k_\mu^{(n)})_{n=0}^\infty$ are called the correlation
functions of the measure $\mu$.

In what follows, we will assume that $k_\mu$ satisfies the Ruelle
bound, i.e., there exists a constant $\xi>0$ such that
\begin{equation}\label{RB} k_\mu(\eta)\le
\xi^{|\eta|}\quad \text{for all } \eta\in\Gamma_0.\end{equation}
Using \eqref{RB}, one, in particular, gets that all local moments of
$\mu$ are finite:
\begin{eqnarray}\label{eq19}
\int_\Gamma |\gamma_\Lambda|^n \,\mu(d\gamma) < \infty, \quad
 n \in {\mathbb N}, \ \Lambda \in {\mathcal
O}_{\mathrm c}({\mathbb R}^{d}).
\end{eqnarray}

We will also use the following lemma.

\begin{lemma}\label{iutgu}
Let $f:\R^d\to\R$ be a measurable function which is bounded  outside a set $\Lambda\in\mathcal O_{\mathrm c}(\R^d)$ and such that $e^f-1\in L^1(\R^d,dx)$. Let also $g,\,g_1,\,g_2:\R^d\to\R$ be such that $e^fg,\,e^fg_1,\,e^fg_2\in L^1(\R^d,dx)$.
Define functions $G_1$, $G_2$, $G_3$ on $\Gamma_0$
by
\begin{align*}
G_1&=\big((e^f-1)^{\otimes n}\big)_{n=0}^\infty,\\
G_2&=\big(n(e^f-1)^{\otimes(n-1)}\odot (e^fg)\big)_{n=0}^\infty,\\
G_3&=\big(n(n-1)(e^f-1)^{\otimes(n-2)}\odot (e^fg_1)\odot(e^fg_2)\big)_{n=0}^\infty,
\end{align*}
where $\odot$ denotes  symmetric tensor product\rom.
Then, $G_1,G_2,G_3\in L^1(\Gamma_0,
\rho_\mu)$ and
\begin{align}
(KG_1)(\gamma)&=e^{\langle f,\gamma\rangle},\notag\\
(KG_2)(\gamma)&=e^{\langle f,\gamma\rangle}
\langle g,\gamma\rangle,\notag\\
(KG_2)(\gamma)&=e^{\langle f,\gamma\rangle}
\sum_{x_1\in\gamma}\sum_{x_2\in\gamma,\, x_2\ne x_1} g_1(x_1)g_2(x_2),\label{hgyu}
\end{align}
for $\mu$-a.e.\ $\gamma\in\Gamma$, and so $KG_1,KG_2,KG_3\in L^1(\Gamma,\mu)$.

\end{lemma}

\noindent {\it Proof}.
Using the Ruelle bound, we clearly have that $G_1,G_2,G_3\in L^1(\Gamma_0,
\rho_\mu)$. Hence, by Proposition~\ref{dftrert}, we get $KG_1,KG_2,KG_3\in L^1(\Gamma,\mu)$.

Since $f$ is bounded on $\Lambda^c$ and
$e^f-1\in L^1(\R^d,dx)$, we have $f\in L^1(\Lambda^c,dx)$. Therefore,
again using the Ruelle bound,
we get: $\la |f|,\gamma_{\Lambda^c}\ra\in L^1(\Gamma,\mu)$. Hence,
$\la |f|,\gamma\ra<\infty$ for $\mu$-a.e.\ $\gamma\in\Gamma$.
Furthermore, we have $g,g_1,g_2\in L^1(\Lambda^c,dx)$, and so the functions $\la |g|,\gamma\ra$
and $\sum_{x_1\in\gamma}\sum_{x_2\in\gamma,\, x_2\ne x_1}|g_1(x_1) g_2(x_2)|\, $ are finite for $\mu$-a.e.\  $\gamma\in\Gamma$.
Thus the functions on the right hand side of formulas \eqref{hgyu} are well-defined and finite for $\mu$-a.e.\ $\gamma\in\Gamma$.

Next, assume that $f,\,g,\,g_1,\,g_2$ have compact support. Then,
equalities  \eqref{hgyu} follow by a straightforward calculation. The general case follows by approximation.\quad $\square$

We introduce a $\star$-convolution of two functions on $\Gamma_0$, so that
$$ (K(G_1\star G_2))(\gamma)=(KG_1)(\gamma)(KG_2)(\gamma) $$
(cf.\ \cite{KK}). Then, we have:
\begin{equation}\label{star} (G_1\star G_2)(\eta)=\sum_{(\eta_1,\eta_2,\eta_3)\in\mathcal P_3(\eta)}G_1(\eta_1\cup\eta_2)G_2(\eta_2\cup\eta_3),\end{equation} where $\mathcal P_3(\eta)$ is the set of all ordered partitions of $\eta$ into three parts.

For each  $\Lambda\subset \R^d$, we denote $$\Gamma_\Lambda
:=\{\gamma\in\Gamma:\, \gamma\subset \Lambda \}.
$$ A measurable function $F:\Gamma\to\R$ is called local if there exists
 $\Lambda\in\mathcal O_{\mathrm c}(\R^d)$ such that
$$ F(\gamma)=F(\gamma_\Lambda)\quad \text{for all }\gamma\in\Gamma. $$
For such a function $F$, the pre-image of $F$ under $K$ is given by
\begin{equation}\label{ghfty} (K^{-1}F)(\eta)=
\chi_{\Gamma_\Lambda}(\eta)\sum_{\xi\subset\eta}
(-1)^{|\eta\setminus\xi|} F(\xi),\end{equation} see e.g.\ \cite{KK}.

We will also need the space $\ddot\Gamma:=\ddot\Gamma_{\R^d}$ which
consists of all multiple configurations in  $\R^d$. So,
$\ddot\Gamma$ is the set of all Radon $\Z_+\cup\{\infty\}$-valued
measures on $\R^d$. In particular, $\Gamma\subset\ddot\Gamma$.
Analogously to the case of $\Gamma$, we define the vague topology on
$\ddot\Gamma$ and the corresponding Borel $\sigma$-algebra $\mathcal
B(\ddot\Gamma)$. For each  $\Lambda\subset \R^d$, we denote
$$\ddot\Gamma_\Lambda :=\{\gamma\in\ddot\Gamma:\,
\operatorname{supp}(\gamma)\subset \Lambda \}.
$$
Also, by analogy, we will say that a measurable function
$F:\ddot\Gamma\to\R$ is local if there exists $\Lambda\in\mathcal
O_{\mathrm{c}}(\R^d)$ such that \begin{equation}\label{fgdgdsdr}
F(\gamma)=F(\gamma_\Lambda)\quad \text{for all
}\gamma\in\ddot\Gamma, \end{equation} where
$\gamma_\Lambda(dx):=\chi_\Lambda(x)\,\gamma(dx)$.

\section{Gibbs measures on configuration spaces}
\label{wqwwwq}

A pair potential (without hard core) is a
Borel  measurable function $\phi\colon \R^d\to \R\cup\{+\infty\}$ such that $\phi(-x)=\phi(x)\in\R$ for all $x\in\R^d\setminus\{0\}$. For $\gamma\in\Gamma$ and $x\in \R^d\setminus\gamma $, we define the relative
energy of interaction between a particle at $x$ and the
configuration $\gamma$ as follows:
\begin{equation}\label{ree}
E(x,\gamma):=\left\{
\begin{aligned}
&\sum_{y\in\gamma}\phi(x-y),&& \text{if } \sum_{y\in\gamma}|\phi(x-y)|
< +\infty ,\\
& +\infty, &&\text{otherwise.}
\end{aligned}
\right.
\end{equation}

A probability measure $\mu$ on $(\Gamma,\mathcal{B}(\Gamma))$ is
called a (grand canonical) Gibbs measure corresponding to
the pair potential $\phi$ and activity $z>0$ if it satisfies
the Georgii--Nguyen--Zessin identity \cite[Theorem 2]{NZ}:
\begin{equation}\label{GNZ}
\int_\Gamma\mu(d\gamma)\int_{\R^d}\gamma(dx) F(\gamma,x)
=\int_\Gamma \mu(d\gamma)\int_{\R^d} z\,dx\, \exp\left[-E(x,\gamma)\right]
F(\gamma\cup x,x)
\end{equation}
for any measurable function $F:\Gamma\times\R^d\to[0,+\infty]$.
We denote the set of all such measures $\mu$ by $\mathcal G (z,\phi)$.

Note that, by virtue of \eqref{ree} and by applying
\eqref{GNZ} twice,  we get, for any measurable function $U:\Gamma\times(\R^d)^2\to[0,+\infty]$,
\begin{multline}\label{jhgftyfwewe}
\int_\Gamma\mu(d\gamma)\sum_{x_1\in\gamma}
\sum_{x_2\in\gamma,\, x_2\ne x_1}U(\gamma,x_1,x_2)
=\int_\Gamma\mu(d\gamma)\int_{\R^d}z\,dx_1
\int_{\R^d}z\,dx_2 \\ \times\exp[-E(x_1,\gamma)-E(x_2,\gamma)-\phi(x_1-x_2)] U(\gamma\cup x_1\cup x_2,x_1,x_2).
\end{multline}

Let us now describe the class of Gibbs measures of Ruelle type  \cite{Ru70}.  We will first formulate conditions on the interaction.

For every $r=(r^1,\dots,r^d)\in\Z^d$, we define the cube
\begin{equation}\label{drtse}Q_r:=\left\{\, x\in\R^d\mid r^i-\frac 12\le x^i<r^i+\frac12
\,\right\}.\end{equation} These cubes form a partition of $\R^d$.
For any $\gamma\in\Gamma$, we set
\begin{equation}\label{drtsghgh}\gamma_r:=\gamma_{Q_r},\quad r\in\Z^d.\end{equation}

\begin{description}

\item[(SS)] ({\it Superstability})
 There exist $A>0$ and  $B\ge0$ such that, for each  $\gamma\in\Gamma_0$,
$$\sum_{\{x,y\}\subset\gamma} \phi(x-y)\ge\sum_{r\in\Z^d}\big(A|\gamma_r|
^2-B|\gamma_r|\big).$$

\end{description}

Notice that the superstability condition automatically implies that the
potential $\phi$ is semi-bounded from below.

\begin{description}

\item[(LR)] ({\it Lower regularity}) There exists a decreasing positive
function $a\colon\N\to{\Bbb R}_+$ such that
$$\sum_{r\in\Z^d}a(\|r\|)<\infty$$ and for any
$\Lambda',\Lambda''$ which are finite unions of cubes $Q_r$ and
disjoint, with $\gamma'\in\Gamma_{\Lambda'}$,
$\gamma''\in\Gamma_{\Lambda''}$,
$$\sum_{x\in\gamma',\, y\in\gamma''}\phi(x-y)\ge-\sum_{r',r''\in\Z^d}a(\|r'-r''\|)
|\gamma_{r'}'|\,|\gamma_{r''}''|.$$ Here, $\|\cdot\|$ denotes the
maximum norm on $\R^d$.

\item[(I)] ({\it Integrability}) We have
$$\int_{\R^d}|e^{-\phi(x)}-1|\, dx<+\infty.$$

\end{description}

For $N\in\N$, let $\Lambda_N$ be the cube with side length $2N-1$ centered at the origin in $\R^d$, $\Lambda_N$ is then a union of $(2N-1)^d$ unit cubes of the form $Q_r$. 

A probability measure $\mu$ on $(\Gamma,{\cal B}(\Gamma))$ is
called tempered if $\mu$ is supported by
$S_\infty{:=}\bigcup_{n=1}^\infty S_n$, where $$S_n:=\left\{\,
\gamma\in\Gamma\mid \forall N\in\N\ \sum_{r\in\Lambda_N\cap\Z^d}
|\gamma_r|^2\le n^2|\Lambda_N \cap\Z^d| \,\right\}.$$ By ${\cal
G}^t(z,\phi)\subset{\cal G}(z,\phi)$ we denote the set of all
tempered grand canonical Gibbs measures.

\begin{theorem}[\cite{Ru70}]\label{ioh}
 Let \rom{(SS), (I),} and  \rom{(LR)} hold. Then the
set ${\cal G}^t(z,\phi)$ is non-empty for each $z>0$. Furthermore, each  $\mu\in {\cal G}^t(z,\phi)$ has correlation functions
 which satisfy the following bound\rom: there exists 
 $\xi,\psi>0$ such that  
 \begin{equation}\label{uigz}
k_\mu(\eta)\le \xi^{|\eta|}\exp\bigg[-\psi\sum_{r\in\Z^d} |\eta_r|^2
\bigg],\quad \text{\rom{for all} }\eta\in\Gamma_0.
\end{equation}

\end{theorem}

Note that  the estimate \eqref{uigz} is evidently stronger than  the Ruelle bound \eqref{RB}.

In what follows, we will keep a Gibbs measure $\mu\in\mathcal G^t(z,\phi)$  as in Theorem~\ref{ioh} fixed, and we
will additionally assume that there exists  $\Theta\in\mathcal O_{\mathrm c}(\R^d)$ such that
\begin{equation}\label{jhf66} \sup_{x\in\Theta^c}\phi(x)<\infty.\end{equation}
Since $\phi$ is bounded from below, (I) is now equivalent to the condition
$\phi\in L^1(\Theta^c,dx)$. Furthermore, by \cite[Lemma~3.1]{KLRDiffusions}, the relative energy
$E(x,\gamma)$  is finite for $dx\otimes \mu$-a.e.\ $(x,\gamma)\in\R^d\times \Gamma$, as well as
$E(x,\gamma\setminus x)$ is finite for
$\mu$-a.e.\ $\gamma\in\Gamma$ and for all $x\in\gamma$.

\section{Kawasaki dynamics}\label{tyedtwe}

We introduce the set $\mathcal FC_{\mathrm b}(C_0(\R^d),\Gamma)$
of all functions of the form $$ \Gamma\ni\gamma\mapsto F(\gamma)=g(\la \varphi_1,\gamma\ra,\dots,
\la\varphi_N,\gamma\ra),$$ where $N\in\N$, $\varphi_1,\dots,\varphi_N\in C_0(\R^d)$, and $g\in C_{\mathrm b}(\R^N)$,
where $ C_{\mathrm b}(\R^N)$ denotes the set of all continuous bounded functions on $\R^N$.
For each function $F:\Gamma\to\R$, $\gamma\in\Gamma$, and $x,y\in\R^d$, we denote
$$ (D^{-+}_{xy}F)(\gamma):=F(\gamma\setminus x\cup y)-F(\gamma).$$

We fix any  $a:\R^d\to[0,\infty)$ which is bounded and such that
$a\in L^1(\R^d,dx)$ and $a(-x)=a(x)$ for all $x\in\R^d$.
We   define a bilinear form
\begin{multline*}
\mathcal E(F,G):=\frac12\int_\Gamma \mu(d\gamma)\int_{\R^d}\gamma(dx)\int_{\R^d}dy \, a(x-y)\\
\times \exp[(1/2)E(x,\gamma\setminus x)-(1/2)E(y,\gamma\setminus x)] (D_{x,y}^{-+}F)(\gamma)
(D_{x,y}^{-+}G)(\gamma),\end{multline*} where $F,G\in \mathcal FC_{\mathrm b}(C_0(\R^d),\Gamma)$.

The following theorem was proved in \cite{KLR}.

\begin{theorem}\label{hvjfgh}
\rom{(i)}  The bilinear form $(\mathcal E,\mathcal FC_{\mathrm b}(C_0(\R^d),\Gamma))$
 is closable on $L^2(\Gamma,\mu)$ and its closure will be denoted by $(\mathcal E,D
 (\mathcal E))$.

\rom{(ii)} There exists a conservative
Hunt process
\[
\mathbf{M}=\big(\mathbf
{\Omega}, \ \mathbf{F}, \
(\mathbf{F}_t)_{t\geq 0}, \ (\mathbf{\Theta}_t)_{t\geq 0}, \ (\mathbf{X}(t))_{t\geq 0}, \
(\mathbf{P}_\gamma)_{\gamma\in\Gamma}\big)
\]
on $\Gamma$ \rom(see e.g.  \rom{\cite[p.~92]{MaRo})} which is properly associated
with $(\mathcal E,D
 (\mathcal E))$, i.e., for all \rom($\mu$-versions
of\,\rom) $F\in L^2(\Gamma,\mu)$ and all $t>0$ the function
\[
\Gamma\ni\gamma\mapsto(p_tF)(\gamma)
:=\int_{\mathbf{\Omega}} F(\mathbf{X}(t))d\mathbf{P}_\gamma
\]
is an $\mathcal E$-quasi-continuous version of
$\exp\left[-tH\right]F$, where $(H,D(H))$ is the generator of
$(\mathcal E,D
 (\mathcal E))$.  
 In particular, $\mathbf{M}$ has $\mu$ as  invariant
measure.
 $\mathbf{M}$ is up to
$\mu$-equivalence unique \rom(cf.\ \rom{\cite[Chap.\ IV, Sect.\ 6]{MaRo})}.

 \rom{(iii)} We have  $\mathcal FC_{\mathrm b}(C_0(\R^d),\Gamma)\subset D(H)$
 and for any $F\in \mathcal FC_{\mathrm b}(C_0(\R^d),\Gamma)$,
\begin{multline}\label{jkhuioo} (HF)(\gamma)=-\int_{\R^d}\gamma(dx)\int_{\R^d}dy
\,a(x-y) \\ \times \exp[(1/2)E(x,\gamma\setminus x)-(1/2)E(y,\gamma\setminus x)] (D_{x,y}^{-+}F)(\gamma). \end{multline}
\end{theorem}

We will call the process $\mathbf{M}$ from Theorem~\ref{hvjfgh} the Kawasaki dynamics
(of continuous particles).

\begin{remark}\rom{
In Theorem~\ref{hvjfgh} (ii),   $\mathbf{M}$ can be taken canonical,
i.e., $\mathbf{\Omega}$ is the set\linebreak
$D(\left[0,+\infty\right),\Gamma)$ of all {\it c\'adl\'ag\/} functions
$\omega:\left[0,+\infty\right)\to\Gamma$ (i.e., $\omega$
is right continuous on $\left[0,+\infty\right)$ and has
left limits on $(0,+\infty)$),
$\mathbf{X}(t)(\omega)=\omega(t)$, $t\geq0$,
$\omega\in\mathbf{\Omega}$, $(\mathbf{F}_t)_{t\geq0}$
together with $\mathbf{F}$ is the correponding minimum
completed admissible family (cf.\ \cite[Section 4.1]{Fu80}) and $\mathbf{\Theta}_t$, $t\geq0$, are the corresponding
natural time shifts.
}\end{remark}




\section{Gradient stochastic dynamics}\label{ytrtyed}

We denote by $\mathfrak D$ the set of all local
functions $F$ on $\ddot\Gamma$ which satisfy the following
assumptions:
\begin{description}

\item[\rom{(i)}] For each fixed $\gamma\in\ddot\Gamma$, the function
$$\R^d\ni x\mapsto F(\gamma+\varepsilon_x)$$
is twice  continuously differentiable. 

\item[\rom{(ii)}] Let $\Lambda$ be the minimal subset of $\R^d$ which is a finite union of $Q_r$ cubes, and such that \eqref{fgdgdsdr} holds for this set $\Lambda$. Then there exist $\zeta>0$, $\tau\ge0$,
$\sigma\ge0$, and $0<p<1$ (depending on $F$) such that, for each 
$\gamma\in\ddot\Gamma_\Lambda$ and each 
$x\in\Lambda$,
\begin{equation}
\label{hguzffdrtdt} |F(\gamma)|\vee \|\nabla_x F(\gamma+\varepsilon_
x)\|\vee \|\nabla^2_x F(\gamma+\varepsilon_x)\|
\le
\zeta^{|\gamma|}\exp\bigg[\,\tau\bigg(1+\sigma
\sum_{r\in\Z^d}|\gamma_r|^2\bigg)^p\,\bigg]. \end{equation}

\end{description}

Note that $\mathfrak D$, in particular, includes
all local functions on $\ddot\Gamma$ which satisfy (i) and for which the left hand side of \eqref{hguzffdrtdt} is bounded, as a function of $\gamma\in\ddot\Gamma$ and $x\in\R^d$.

We also introduce the set $\mathcal FC^2_{\mathrm b}(C^2_0(\R^d),\ddot\Gamma)$
of all functions of the form
$$ \ddot\Gamma\ni\gamma\mapsto F(\gamma)=g(\la \varphi_1,\gamma\ra,\dots,
\la\varphi_N,\gamma\ra),$$ where $N\in\N$, $\varphi_1,\dots,\varphi_N\in C^2_0(\R^d)$, and $g\in C^2_{\mathrm b}(\R^N)$.
Here and below, $C_0^k(\R^d)$ and $C^k_{\mathrm b}(\R^N)$, $k\in\N$,  denote the space of all $k$ times continuously differentiable functions on $\R^d$ with compact support, respectively the space of all bounded, $k$ times continuously differentiable functions on $\R^N$ with bounded derivatives.
We evidently have the inclusion
$$ \mathcal FC^2_{\mathrm b}(C^2_0(\R^d),\ddot\Gamma)\subset \mathfrak D,$$
and therefore the set $\mathfrak D$ is dense in $L^2(\Gamma,\mu)$.
(We have included functions from $\mathfrak D$ into $L^2(\Gamma,\mu)$ by taking their restriction to  $\Gamma$.)

In what follows, we will use the following 

\begin{lemma}\label{esse} 
Let $\Lambda\subset\R^d$ be a finite union of $Q_r$ cubes and let $\zeta>0$, $\tau\ge0$, $\sigma\ge0$, and $0<p<1$. Define
$$ U(\gamma):=\zeta^{|\gamma_\Lambda|}\exp\bigg[\tau\bigg(1+\sigma\sum_{r\in\Z^d\cap\Lambda}|\gamma_r|^2\bigg)^p\bigg],\qquad\gamma\in\Gamma.$$ Then, for each $\eta\in\Gamma_0$, 
$$ |(K^{-1}U)(\eta)|\le \chi_{\Gamma_\Lambda}(\eta)(2\zeta)^{|\eta|}\exp\bigg[\tau\bigg(1+\sigma\sum_{r\in\Z^d}|\eta_r|^2\bigg)^p\bigg].$$
\end{lemma}

\noindent{\it Proof.} The lemma follows from \eqref{ghfty} if we take into account that the sum in  \eqref{ghfty} has exactly $2^{|\eta|}$ terms.\vspace{2mm} \quad $\square$ 


We fix any $c>0$ and define a bilinear form
\begin{equation}\label{dytt}
\mathcal E^{(\mathrm{dif})}(F,G):=\frac{c}2\int_\Gamma\mu(d\gamma)\int_{\R^d}z\, dx\,
\la\nabla_x F(\gamma+\varepsilon_x),\nabla_x G(\gamma+\varepsilon_x)
\ra\exp\big[-E(x,\gamma)\big],\end{equation}
where $F,G\in \mathfrak D $, and we denoted by $\la\cdot,\cdot\ra$
the scalar product in $\R^d$. 
Using the Cauchy--Schawrz inequality, Theorem~\ref{ioh}, \eqref{hguzffdrtdt}, and Lemma~\ref{esse},  the integral on the right hand side of \eqref{dytt} is well defined and finite.

For a function $F:\Gamma\to\R$,  a fixed $\gamma\in\Gamma$ and $x\in\gamma$, we denote
\begin{equation}\label{fytf} \nabla_x F(\gamma):=\nabla_y F(\gamma-\varepsilon_x+\varepsilon_y)\big|_{y=x},\end{equation}
provided the  gradient on the right hand side of \eqref{fytf} exists at point $x$. Then, by \eqref{GNZ}, we also have 
\begin{equation}\label{hjuuf}
\mathcal E^{(\mathrm{dif})}(F,G)=\frac{c}2\int_\Gamma\mu(d\gamma)\int_{\R^d}\gamma(dx)
\la\nabla_x F(\gamma),\nabla_x G(\gamma)\ra,
\end{equation}
for  $F,G\in \mathfrak D $.

The following theorem follows from  
\cite{AKR4,MR98,RS98}, see also \cite{Osa96,Yos96}.

\begin{theorem}\label{jhvf} Assume that   $\phi$ is  differentiable
on $\R^d\setminus\{0\}$, $e^{-\phi}$ is differentiable on $\R^d$, and we have
$$ \|\nabla\phi\|\in L^1(\R^d,e^{-\phi(x)}\,dx)\cap  L^2(\R^d,e^{-\phi(x)}\,dx).$$

Then\rom:

\begin{description}

\item[\rom{(i)}] The bilinear form $(\mathcal E^{(\mathrm{dif})},  \mathfrak D)$
is closable on $L^2(\Gamma,\mu)$ and its closure will be denoted by $(\mathcal E^{(\mathrm{dif})},D
 (\mathcal E^{(\mathrm{dif})}))$.

 \item[\rom{(ii)}] Denote by $(H^{(\mathrm{dif})},D(H^{(\mathrm{dif})}))$
 the generator of  $(\mathcal E^{(\mathrm{dif})},D
 (\mathcal E^{(\mathrm{dif})}))$. Then \linebreak 
 $\mathfrak D\subset D
 (H^{(\mathrm{dif})})$ and for each $F\in \mathfrak D$,
 \begin{equation}\label{hjfgytu}
( H^{(\mathrm{dif})}F)(\gamma)=\frac{c}2\int_{\R^d}\gamma(dx)
\bigg(- \Delta_xF(\gamma)+\sum_{u\in\gamma\setminus x}\la
\nabla_x F(\gamma),\nabla\phi(x-u)\ra\bigg).
  \end{equation}
  Here, $\Delta_xF(\gamma):=\Delta_uF(\gamma\setminus x\cup u)\big|_{u=x}$.

\item[\rom{(iii)}] There exists a conservative
diffusion process
\begin{equation*}
\mathbf{M}^{(\mathrm{dif})}=\big(\mathbf{\Omega}^{(\mathrm{dif})}, \ \mathbf{F}^{(\mathrm{dif})}, \
(\mathbf{F}^{(\mathrm{dif})}_t)_{t\geq 0}, \ (\mathbf{\Theta}^{(\mathrm{dif})}_t)_{t\geq 0}, \ (\mathbf{X}^{(\mathrm{dif})}(t))_{t\geq 0}, \
(\mathbf{P}^{(\mathrm{dif})}_\gamma)_{\gamma\in\ddot\Gamma}\big)
\end{equation*}
on $\ddot\Gamma$ \rom(see e.g.  \rom{\cite[p.~92]{MaRo})} which is properly associated
with $(\mathcal E^{(\mathrm{dif})},D
 (\mathcal E^{(\mathrm{dif})}))$. 
 In particular, $\mathbf{M}^{(\mathrm{dif})}$ has $\mu$ as  invariant
measure.  The $\mathbf{M}^{(\mathrm{dif})}$ is up to
$\mu$-equivalence unique.

\item[\rom{(iv)}] In the case $d\ge2$, the set
$\ddot\Gamma\setminus\Gamma$ is $\mathcal E^{(\mathrm{dif})}$-exceptional, so that
$\ddot\Gamma$ may be replaced with $\Gamma$ in \rom{(iii)}.
\end{description}

\end{theorem}

\begin{remark} \rom{ Note that, even when $d=1$, the finite-dimensional
distributions of the process $\mathbf {M}^{(\mathrm{dif})}$ are
concentrated on the Cartesian powers of the space $\Gamma$.}
\end{remark}

\begin{remark}\rom{
Note that the initial domain $\mathfrak D$ of the bilinear form $\mathcal E^{(\mathrm{dif})}$
is bigger than the domains of the corresponding bilinear forms in \cite{AKR4,MR98,RS98,Yos96}. 
It is, generally speaking, an open problem whether all these forms coincide after being closed (compare with Remark~4.14 in \cite{MR98}).
Note also that, even in the case of a bilinear form  $\mathcal E^{(\mathrm{dif})}$ with a smaller domain, the convergence result of Theorem~\ref{uyfytdy} below will still be true, however for a smaller set of functions $F$.   

}\end{remark}

\section{Main results}\label{yfyyty}

Let us consider the Kawasaki dynamics
 $\mathbf{M}$ from Theorem~\ref{hvjfgh}.
 We will assume that $a(x)=\tilde a(|x|)$ for all $x\in\R^d$, where $\tilde a:[0,\infty)\to [0,\infty)$.
 We now perform the following scaling of
 this dynamics. For each $\epsilon>0$, instead of the function $a$,
 we use the function $a_\epsilon$ given by \eqref{rtsra}.
   In the obtained dynamics, we also scale  time,
multiplying it  by $\epsilon^{-2}$. Thus, we obtain a  Kawasaki
dynamics
 $\mathbf{M}^{(\epsilon)}$, which is exactly the Hunt process from Theorem~\ref{hvjfgh} corresponding to the function $\epsilon^{-2}a_\epsilon$. We denote by $(H^{(\epsilon)},D(H^{(\epsilon)}))$ the generator of this dynamics.  Completely analogously to the proof of \cite[Lemma~4.1]{KL}, we conclude that, for each $\epsilon>0$,   $\mathfrak D \subset D(H^{(\epsilon)})$ and, for each $F\in \mathfrak D$,
\begin{multline} (H^{(\epsilon)}F)(\gamma)=-\epsilon^{-d-2}\int_{\R^d}\gamma(dx)\int_{\R^d}dy
\,a((x-y)/\epsilon) \\  \times \exp[(1/2)E(x,\gamma\setminus x)-(1/2)E(y,\gamma\setminus x)] (D_{x,y}^{-+}F)(\gamma).\label{ytdrtsd}\end{multline}

\begin{theorem}\label{uyfytdy} 
Let the conditions of   Theorem~\rom{\ref{jhvf}} be satisfied. Furthermore,  assume that the
following conditions are satisfied\rom:

\begin{description}

\item{\rom{a)}} The function $a$ has compact support\rom.

\item{\rom{b)}} We have $e^{-\phi/2}\in C_{\mathrm b}^1(\R^d)$\rom.

\item{\rom{c)}} For each $\delta>0$, set 
$$g_\delta(x):=\sup_{y\in B(x;\delta)}e^{-\phi(y)/2}\|\nabla\phi(y)\|,\qquad x\in\R^d.$$ Here, $B(x;\delta)$ denotes the closed  ball in $\R^d$ centered at $x$ and of radius $\delta$. 
Then, 
then there exists $\delta>0$ such that $g_\delta\in L^1(\R^d,dx)$\rom.

\item{\rom{d)}} There exists $\Lambda\in\mathcal O_{\mathrm c}(\R^d)$ such that  $$ e^{-\phi/2}\|\nabla\phi\|\in L^1(\Lambda,dx)$$
and the function $e^{-\phi/2}\|\nabla\phi\|$ is bounded
on $\Lambda^c$.
\end{description}

Let \begin{equation}\label{dtse} c:=\int_{\R^d}a(x)(x^1)^2\,dx
\end{equation} and let
$(H^{(\mathrm{dif})},D(H^{(\mathrm{dif})}))$
correspond to the above choice of the constant $c$ \rom(see \eqref{hjfgytu}\rom)\rom. Then, for each
$F\in \mathfrak D$, we have\rom:
\begin{equation}\label{ghdcrtrt} H^{(\epsilon)} F\to H^{(\mathrm{dif})}F\quad\text{\rom{in} $L^2(\Gamma,\mu)$ \rom{as} $\epsilon\to0$}.\end{equation}

\end{theorem}

\begin{remark}\rom{ Note that condition a) heuristically means that, in the initial Kawasaki dynamics, there is a finite maximal length of jumps of particles.   Notice also
that condition c) of Theorem \ref{uyfytdy} is slightly stronger than the condition
$e^{-\phi/2}\|\nabla\phi\|\in L^1 (\R^d,dx)$.
}\end{remark}

Next, we take the canonical realizations of the
processes $\mathbf{M}^{(\epsilon)}$, $\epsilon>0$, and $\mathbf{M}^{( \mathrm{dif})}$ and define
stochastic processes
$\mathbf Y^{(\epsilon)}=(\mathbf Y^{(\epsilon)}_t )_{t\ge0}$
and $\mathbf{Y}^{(\mathrm{dif})}=(\mathbf Y^{(\mathrm{dif})}_t)_{t\ge0}$ whose law is the probability measure on
$D([0,+\infty),\Gamma)$, respectively $C([0,+\infty),\Gamma)$ (replace $\Gamma$ with $\ddot\Gamma$ if $d=1$), given by
$$ \mathbf Q^{(\epsilon)}:=\int_\Gamma \mathbf P^{(\epsilon)}_\gamma\,\mu(d\gamma),$$
respectively
$$ \mathbf Q^{(\mathrm{dif})}:=\int_\Gamma \mathbf P^{(\mathrm{dif})}_\gamma\,\mu(d\gamma).$$

\begin{corollary}\label{iuguyfd}
Assume that the conditions of Theorem~\rom{\ref{uyfytdy}} are satisfied. Assume additionally that $\mathfrak D$ is a core for $(H^{(\mathrm{dif})},D(H^{(\mathrm{dif})}))$\rom.
Then, as $\epsilon\to0$, the finite-dimensional
distributions of the process $\mathbf{M}^{(\epsilon)}$ weakly
converge to the finite-dimensional distributions of the process
$\mathbf{M}^{(\mathrm{dif})}$ with $c$ given by \eqref{dtse}.

\end{corollary}

Following \cite{CPY}, we will now introduce additional conditions on the
potential $\phi$.

Let $\alpha:[0,\infty)\to\R$ be any monotonic, increasing, and concave function such that:

\begin{description}

\item[\rom{(i)}] $\alpha(0)\ge1$ and $\alpha(\lambda)\to\infty$ as $\lambda\to\infty$.

\item[\rom{(ii)}] $\alpha'(\lambda)\le[1/(1+\lambda)]\alpha(\lambda)$ for $\lambda\ge0$, and there exists a constant $c>0$ such that $\alpha''(\lambda)\ge -c[1/(1+\lambda)]$.

\end{description}

For example, let $l(\lambda):=\log(1+\lambda)$, $\lambda\ge0$.
Then, for any $n\in\N$, the function
$\alpha(\lambda):=1+l\underbrace{\circ\dots\circ}_{\text{$n$ times}}l(\lambda)$
 satisfies the above conditions.

So, in what follows we will assume:

\begin{description}
\item[(A)] We have $\phi\in C_{\mathrm b}^3(\R^d)$, and there exist a
constant $c_0$ and a function $\alpha$ that satisfies the conditions
(i) and (ii) above, such that,  for all $x\in\R^d$,
$$ \|\nabla\phi(x)\|+\|\nabla^2\phi(x)\|+\|\nabla^3\phi(x)\|\le
\exp[-c_0\log(1+|x|^2)\alpha(1+|x|^2)].$$
\end{description}

It was proved in \cite{CPY} that, under  condition (A), the set $\mathfrak D$ is a core for the operator $(H^{(\mathrm{dif})},D(H^{(\mathrm{dif})}))$. (In fact, Choi {\it et al.}\ \cite{CPY} found a core for the operator $(H^{(\mathrm{dif})},D(H^{(\mathrm{dif})}))$ which is, as can be easily checked, a subset of $\mathfrak D$.) Furthermore, under condition (A), the potential $\phi$ clearly satisfies assumptions of Theorem~\ref{uyfytdy}.
Thus, we get from Corollary~\ref{iuguyfd}:

\begin{corollary}\label{gtfdytsee}  Assume that the function $a$ has compact support, and assume that
condition \rom{(A)} is satisfied. Then, as $\epsilon\to0$, the finite-dimensional
distributions of the process $\mathbf{M}^{(\epsilon)}$ weakly
converge to the finite-dimensional distributions of the process
$\mathbf{M}^{(\mathrm{dif})}$ with $c$ given by \eqref{dtse}.
\end{corollary}

\begin{remark}
\rom{ Let us briefly explain a generalization of 
Theorem \ref{uyfytdy}. Let us fix a parameter $s\in[0,1]$. (Note that the results of this paper  will correspond to the choice of parameter  $s=1/2$.) 
By \cite{KLR}, there exists a conservative Hunt process on $\Gamma$ (a Kawasaki dynamics) whose $L^2(\Gamma,\mu)$-generator $(H_s,D(H_s))$ is the Friedrichs extension of the operator $(H_s,\mathfrak D)$ given by   
\begin{multline*} (H_{s}F)(\gamma)=-\int_{\R^d}\gamma(dx)\int_{\R^d}dy
\,a(x-y) \\ \times \exp[(1-s)E(x,\gamma\setminus x)-sE(y,\gamma\setminus x)] (D_{x,y}^{-+}F)(\gamma),\qquad F\in\mathfrak D \end{multline*}
(in the case $s<1/2$, the potential $\phi$ must satisfy an additional assumption which reduces the ``strength of singularity'' at zero).

Next, for each $c>0$, it can be shown that, under some conditions on $\phi$ which are analogous to the conditions of Theorem~\ref{jhvf}, 
there exists a conservative diffusion process on $\Gamma$ (respectively, on $\ddot\Gamma$ if $d=1$) whose $L^2(\Gamma,\mu)$-generator 
$(H_s^{(\mathrm{dif})},D(H_s^{(\mathrm{dif})}))$ is the Friedrichs extension of the operator $(H_s^{(\mathrm{dif})},\mathfrak D)$ given by   
 \begin{multline*}\label{hjfgytu}
( H_s^{(\mathrm{dif})}F)(\gamma)=c\int_{\R^d}\gamma(dx)
\bigg(-\frac12\, \Delta_xF(\gamma)+\sum_{u\in\gamma\setminus x}\la
\nabla_x F(\gamma),s\nabla\phi(x-u)\ra\bigg)\\
\times\exp\bigg[ (-2s+1)\sum_{v\in\gamma\setminus x}\phi(x-v)\bigg],\qquad F\in\mathfrak D.
  \end{multline*} Note that, in the case $s=0$, such a diffusive dynamics has been considered in 
\cite{KLRDiffusions}. 

Then, under the same scaling of the Kawasaki dynamics 
$$ a(\cdot)\mapsto \epsilon^{-d-2}a(\cdot/\epsilon),$$
 and with the same choice  of the constant $c$,  \eqref{dtse}, we get the convergence of the generators on $\mathfrak D$. More precisely, under an appropriate modification of the conditions  
 Theorem~\ref{uyfytdy}, we get 
 for each $F\in\mathfrak D$:
\begin{equation}\label{drdrtdrt} H_s^{(\epsilon)} F\to H_s^{(\mathrm{dif})}F\quad\text{\rom{in} $L^2(\Gamma,\mu)$ \rom{as} $\epsilon\to0$}\end{equation} (we have used the obvious notation $H_s^{(\epsilon)}$). 

As for weak convergence of finite-dimensional distributions of the corresponding equilibrium processes, it will follow from \eqref{drdrtdrt} if $\mathfrak D$ is a core for  $(H_s^{(\mathrm{dif})},D(H_s^{(\mathrm{dif})}))$. However, in the case $s\ne1/2$, no result has yet been proved about a core for this generator.

}\end{remark}

\section{Proofs}\label{jhgftfru}

{\it Proof of  Theorem} \ref{uyfytdy}. Denote the support of the
function $a$ by $\Delta$. By a), the set $\Delta$ is bounded and
hence $r:=\sup_{h\in\Delta}|h|<\infty$. Recall $\delta$ from
condition c) of the theorem. In what follows, we will assume that
$\epsilon\in(0,\delta/r)$. Then
\begin{equation}\label{ghjfdrtdseser}|\epsilon h|<\delta\quad \text{for all }h\in\Delta.\end{equation}

Fix any $F\in\mathfrak D$. By  \eqref{ytdrtsd},
\begin{multline*} (H^{(\epsilon)}F)(\gamma)=-\epsilon^{-2}\int_{\R^d}\gamma(dx)\int_{\Delta}dh
\,a(h) \notag \\  \times \exp[(1/2)E(x,\gamma\setminus x)-(1/2)E(x+\varepsilon h,\gamma\setminus x)] (F(\gamma\setminus x\cup(x+\varepsilon h))-F(\gamma)).
\end{multline*}
Using \eqref{GNZ} and \eqref{jhgftyfwewe}, we have:
\begin{align}
&\int_\Gamma (H^{( \epsilon)}F)^2(\gamma)\,\mu(d\gamma)\notag\\
&\quad =\epsilon^{-4}\int_\Gamma \mu(d\gamma)
\int_{\R^d} z\,dx\int_{\Delta}dh_1\int_{\Delta}dh_2\,
a(h_1)a(h_2)
(F(\gamma\cup(x+\epsilon h_1))-F(\gamma\cup x))\notag\\
&\qquad\times
(F(\gamma\cup(x+\epsilon h_2))-F(\gamma\cup x))\notag\\
&\qquad\times \exp\bigg[
\sum_{u\in\gamma}(-(1/2)\phi(x+\epsilon h_1-u)-(1/2)\phi(x+\epsilon h_2-u))
\bigg]\notag\\
&\qquad+ \epsilon^{-4}\int_\Gamma \mu(d\gamma)
\int_{\R^d} z\,dx_1\int_{\R^d}z\,dx_2\int_{\Delta}dh_1\int_{\Delta}dh_2\, a(h_1)a(h_2)\notag\\
&\qquad\times(F(\gamma\cup(x_1+\epsilon h_1)\cup x_2)-F(\gamma\cup x_1\cup x_2))\notag\\
&\qquad \times(F(\gamma\cup x_1\cup (x_2+\epsilon h_2))-F(\gamma\cup x_1\cup x_2))\notag\\
&\qquad\times\exp\bigg[ 
-(1/2)\phi(x_1+\epsilon h_1-x_2)-(1/2)\phi(x_2+\epsilon h_2-x_1)\notag\\&\qquad
+
\sum_{u\in\gamma}(-(1/2)\phi(x_1-u)-(1/2)\phi(x_2-u)
-(1/2)\phi(x_1+\epsilon h_1-u)\notag\\
&\qquad\text{}-(1/2)\phi(x_2+\epsilon h_2-u))\bigg].
\label{ftydr}\end{align}
Here and below, our calculations are justified by the assumptions of the theorem, the definition of $\mathfrak D$,   Lemma~\ref{iutgu},
\eqref{star}, \eqref{uigz}, and Lemma~\ref{esse}. Hence, by \eqref{ftydr} and Lemma~\ref{iutgu}, we get:
\begin{align}
&\int_\Gamma (H^{( \epsilon)}F)^2(\gamma)\,\mu(d\gamma)\notag\\
&\quad =\epsilon^{-4}
\int_{\R^d} z\,dx\int_{\Delta}dh_1\int_{\Delta}dh_2\,
a(h_1)a(h_2)\int_\Gamma \mu(d\gamma)
(F(\gamma\cup(x+\epsilon h_1))-F(\gamma\cup x))\notag\\
&\qquad\times
(F(\gamma\cup(x+\epsilon h_2))-F(\gamma\cup x))\notag\\
&\qquad\times K\left(\big(\big(e^{
-(1/2)\phi(x+\epsilon h_1-\cdot)-(1/2)\phi(x+\epsilon h_2-\cdot)}-1\big)^{\otimes n}\big)_{n=0}^\infty\right)(\gamma)
\notag\\
&\qquad+ \epsilon^{-4}
\int_{\R^d} z\,dx_1\int_{\R^d}z\,dx_2\int_{\Delta}dh_1\int_{\Delta}dh_2\, a(h_1)a(h_2)\notag\\
&\qquad\times
\exp\big[ 
-(1/2)\phi(x_1+\epsilon h_1-x_2)-(1/2)\phi(x_2+\epsilon h_2-x_1) \big]\notag\\
&\qquad\times\int_\Gamma \mu(d\gamma)(F(\gamma\cup(x_1+\epsilon h_1)\cup x_2)-F(\gamma\cup x_1\cup x_2))\notag\\
&\qquad \times(F(\gamma\cup x_1\cup (x_2+\epsilon h_2))-F(\gamma\cup x_1\cup x_2))\notag\\
&\qquad\times
K\left(\big(
\big(e^{-(1/2)\phi(x_1-\cdot)-(1/2)\phi(x_2-\cdot)
-(1/2)\phi(x_1+\epsilon h_1-\cdot)
-(1/2)\phi(x_2+\epsilon h_2-\cdot)}-1
\big)^{\otimes n}\big)_{n=0}^\infty\right)(\gamma)
.\label{gdtrs}
\end{align}

For each $\gamma\in\Gamma$ and $x,h\in\R^d$,
denote by $y_1(\gamma,x,h)$  a point in the segment $[x,x+h]$ such that
\begin{equation}\label{hgxsddes} F(\gamma\cup(x+h))-F(\gamma\cup x)=\langle
\nabla_x F(\gamma\cup x),h\rangle+
\frac12\langle\nabla^2_yF(\gamma\cup y),h^{\otimes 2}\rangle\big|_{y=y_1(\gamma,x,h)}.\end{equation}
Also, for each $x,h\in\R^d$, we denote
by $y_2(x,h)$ a point in the segment $[x,x+h]$ such that
\begin{equation}\label{bhdrt} e^{-(1/2)\phi(x+h)}=e^{-(1/2)\phi(x)}+e^{-(1/2)\phi(y_2(x,h))}\langle-(1/2)\nabla\phi(y_2(x,h)),h\rangle.\end{equation}
Note that the existence of $y_1(\gamma,x,h)$ and $y_2(x,h)$ follows from the definition  of the $\mathfrak D$ and  assumption b) of the theorem, respectively. Note also that, by \eqref{ghjfdrtdseser},
$$ e^{-(1/2)\phi(y_2(x,\epsilon h))}\|\nabla
\phi(y_2(x,\epsilon h))\|\le g_\delta(x),\quad
x\in\R^d,\ h\in\Delta.$$

Now, by \eqref{gdtrs}, \eqref{hgxsddes}, and \eqref{bhdrt}, we have: \begin{align}
&\int_\Gamma (H^{(\epsilon)}F)^2(\gamma)\,\mu(d\gamma) \notag\\
&\quad =
\int_{\R^d} z\,dx\int_{\Delta}dh_1\int_{\Delta}dh_2\,
a(h_1)a(h_2)\int_\Gamma \mu(d\gamma)
\big(\epsilon^{-2}F_{-2}^{(1)}(\gamma,x, h_1, h_2) \notag \\
&\qquad+\epsilon^{-1}F_{-1}^{(1)}(\gamma,x, h_1, h_2,\epsilon)+F_{0}^{(1)}(\gamma,x, h_1, h_2,\epsilon) \big)
 (KG^{(1)}(\cdot,x, h_1,h_2,\epsilon))(\gamma)\notag\\
&\qquad+
\int_{\R^d}z\,dx_1\int_{\R^d}\,z\,dx_2
\int_{\Delta}dh_1\int_\Delta dh_2\, a(h_1)a(h_2)\notag\\
&\qquad\times\big(
e^{-\phi(x_1-x_2)}+\epsilon u_1(x_1,x_2,h_1,h_2,\epsilon)+\epsilon^2u_2(x_1,x_2,h_1,h_2,\epsilon)\big)\notag\\
&\qquad\times
\int_\Gamma \mu(d\gamma)\big(\epsilon^{-2}F_{-2}^{(2)}(\gamma,x_1,x_2, h_1, h_2) +\epsilon^{-1}F_{-1}^{(2)}(\gamma,x_1,x_2, h_1, h_2,\epsilon)\notag\\
&\qquad+F_{0}^{(2)}(\gamma,x_1,x_2, h_1, h_2,\epsilon) \big)(KG^{(2)}(\cdot,x_1,x_2,h_1,h_2,\epsilon))(\gamma).
\label{hgdtrdswewe}
\end{align}
Here,
\begin{align}
& u_1(x_1,x_2,h_1,h_2,\epsilon):=e^{-(1/2)\phi(x_1-x_2)}\big(e^{-(1/2)\phi(y_2(x_1-x_2,\epsilon h_1))}\langle -(1/2)\nabla\phi(y_2(x_1-x_2,\epsilon h_1)),h_1\rangle\notag\\
&\qquad+e^{-(1/2)\phi(y_2(x_2-x_1,\epsilon h_2))}\langle -(1/2)\nabla\phi(y_2(x_2-x_1,\epsilon h_2)),h_2\rangle\big),\notag\\
& u_2(x_1,x_2,h_1,h_2,\epsilon):=e^{-(1/2)\phi(y_2(x_1-x_2,\epsilon h_1))}e^{-(1/2)\phi(y_2(x_2-x_1,\epsilon h_2))}\notag\\
&\qquad\times\langle -(1/2)\nabla\phi(y_2(x_1-x_2,\epsilon h_1)),h_1\rangle
\langle -(1/2)\nabla\phi(y_2(x_2-x_1,\epsilon h_2)),h_2\rangle,\label{hbguy}
\end{align}
and
\begin{align}
&F^{(1)}_{-2}(\gamma,x, h_1, h_2):= \langle \nabla_x F(\gamma\cup x),h_1\rangle
\la\nabla_x F(\gamma\cup x),h_2\rangle,\notag\\
&F^{(1)}_{-1}(\gamma,x, h_1, h_2,\epsilon):=
\langle \nabla_x F(\gamma\cup x),h_1\rangle(1/2)
\langle \nabla_y^2F(\gamma\cup y),h_2^{\otimes 2}\rangle\big|_{y=y_1(\gamma,x,\epsilon h_2)}\notag\\
&\qquad+ \langle \nabla_x F(\gamma\cup x),h_2\rangle(1/2)
\langle \nabla_y^2F(\gamma\cup y),h_1^{\otimes 2}\rangle\big|_{y=y_1(\gamma,x,\epsilon h_1)},\notag\\
&F^{(1)}_{0}(\gamma,x, h_1, h_2,\epsilon)=(1/4)
\langle \nabla_y^2F(\gamma\cup y),h_1^{\otimes 2}\rangle\big|_{y=y_1(\gamma,x,\epsilon h_1)}\notag\\
&\qquad \times \langle \nabla_y^2F(\gamma\cup y),h_2^{\otimes 2}\rangle\big|_{y=y_1(\gamma,x,\epsilon h_2)},\label{rewa}
\end{align}
and
\begin{align}
&F^{(2)}_{-2}(\gamma,x_1,x_2, h_1, h_2):= \langle \nabla_{x_1} F(\gamma\cup x_1\cup x_2),h_1\rangle
\la\nabla_{x_2} F(\gamma\cup x_1\cup x_2),h_2\rangle,\notag\\
&F^{(2)}_{-1}(\gamma,x_1,x_2, h_1, h_2,\epsilon)\notag\\
&\quad:=
\langle \nabla_{x_1} F(\gamma\cup x_1\cup x_2),h_1\rangle(1/2)
\langle \nabla_y^2F(\gamma\cup x_1\cup y),h_2^{\otimes 2}\rangle\big|_{y=y_1(\gamma\cup x_1,x_2,\epsilon h_2)}\notag\\
&\qquad+ \langle \nabla_{x_2} F(\gamma\cup x_1\cup x_2),h_2\rangle(1/2)
\langle \nabla_y^2F(\gamma\cup y\cup x_2),h_1^{\otimes 2}\rangle\big|_{y=y_1(\gamma\cup x_2,x_1,\epsilon h_1)},\notag\\
&F^{(1)}_{0}(\gamma,x, h_1, h_2,\epsilon):=(1/4)
\langle \nabla_y^2F(\gamma\cup y\cup x_2),h_1^{\otimes 2}\rangle\big|_{y=y_1(\gamma\cup x_2,x_1,\epsilon h_1)}\notag\\
&\qquad \times \langle \nabla_y^2F(\gamma\cup x_1\cup y),h_2^{\otimes 2}\rangle\big|_{y=y_1(\gamma\cup x_1,x_2,\epsilon h_2)},\label{jkhbkigba}
\end{align}
and
\begin{align}
&G^{(1)}(\cdot,x,h_1, h_2,\epsilon)\notag\\
&\quad := \big((e^{-\phi(x-\cdot)}-1+\epsilon g^{(1)}_1(\cdot,x,h_1, h_2,\epsilon)
+\epsilon^2 g^{(1)}_2(\cdot,x,h_1, h_2,\epsilon))^{\otimes n}\big)_{n=0}^{\infty},\notag\\
&G^{(2)}(\cdot,x_1,x_2,h_1, h_2,\epsilon)
:=\big((e^{-\phi(x_1-\cdot)-\phi(x_2-\cdot)}-1\notag\\
&\qquad+\epsilon g^{(2)}_1(\cdot,x_1,x_2,h_1, h_2,\epsilon)
+\epsilon^2 g^{(2)}_2(\cdot,x_1,x_2,h_1, h_2,\epsilon))^{\otimes n}\big)_{n=0}^{\infty},\notag\end{align}
where
\begin{align}
&g^{(1)}_1(\cdot,x,h_1, h_2,\epsilon):=
e^{-(1/2)\phi(x-\cdot)}\big(e^{-(1/2)\phi(y_2(x-\cdot,\epsilon h_1))}\langle -(1/2)\nabla\phi(y_2(x-\cdot,\epsilon h_1)),h_1\rangle\notag\\
&\qquad+e^{-(1/2)\phi(y_2(x-\cdot,\epsilon h_2))}\langle -(1/2)\nabla\phi(y_2(x-\cdot,\epsilon h_2)),h_2\rangle\big),\notag\\
&g^{(1)}_2(\cdot,x,h_1, h_2,\epsilon):=e^{-(1/2)\phi(y_2(x-\cdot,\epsilon h_1))}\langle -(1/2)\nabla\phi(y_2(x-\cdot,\epsilon h_1)),h_1\rangle\notag\\
&\qquad\times
e^{-(1/2)\phi(y_2(x-\cdot,\epsilon h_2))}\langle -(1/2)\nabla\phi(y_2(x-\cdot,\epsilon h_2)),h_2\rangle,\notag\\
&g_1^{(2)}(\cdot,x_1,x_2,h_1, h_2,\epsilon):=
e^{-(1/2)\phi(x_1-\cdot)-\phi(x_2-\cdot)}
e^{-(1/2)\phi(y_2(x_1-\cdot,\epsilon h_1))}\notag\\
&\qquad\times\langle -(1/2)\nabla\phi(y_2(x_1-\cdot,\epsilon h_1)),h_1\rangle\notag\\
&\qquad +e^{-\phi(x_1-\cdot)-(1/2)\phi(x_2-\cdot)}
e^{-(1/2)\phi(y_2(x_2-\cdot,\epsilon h_2))}\langle -(1/2)\nabla\phi(y_2(x_2-\cdot,\epsilon h_2)),h_2\rangle,\notag\\
&g_2^{(2)}(\cdot,x_1,x_2,h_1, h_2,\epsilon):=
e^{-(1/2)\phi(x_1-\cdot)-(1/2)\phi(x_2-\cdot)}
e^{-(1/2)\phi(y_2(x_1-\cdot,\epsilon h_1))}
e^{-(1/2)\phi(y_2(x_2-\cdot,\epsilon h_2))}\notag\\
&\qquad\times \langle -(1/2)\nabla\phi(y_2(x_1-\cdot,\epsilon h_1)),
h_1\rangle
\langle -(1/2)\nabla\phi(y_2(x_2-\cdot,\epsilon h_2)),
h_2\rangle
.\label{hyfyt}
\end{align}

Since $F$ is a local function, so are $F_j^{(i)}$,
$i=1,2$, $j=-2,-1,0$, as functions of $\gamma\in\Gamma$.
Then, by \eqref{hgdtrdswewe}, we get
\begin{align}
&\int_\Gamma (H^{( \epsilon)}F)^2(\gamma)\,\mu(d\gamma) \notag\\
&\quad =
\int_{\R^d} z\,dx\int_{\Delta}dh_1\int_{\Delta}dh_2\,
a(h_1)a(h_2)\notag\\
&\quad\times
\int_{\Gamma_0}\rho_\mu(d\eta)\big(K^{-1}\big(\epsilon^{-2}
F^{(1)}_{-2}(\cdot,x, h_1, h_2)
+\epsilon^{-1} F^{(1)}_{-1}(\cdot,x, h_1, h_2,\epsilon)\notag\\
&\qquad + F^{(1)}_{0}(\cdot,x, h_1, h_2,\epsilon)\big)  \star
G^{(1)}(\cdot,x, h_1,h_2,\epsilon)\big)(\eta)\notag\\
&\qquad+ \int_{\R^d}z\,dx_1\int_{\R^d}\,z\,dx_2
\int_{\Delta}dh_1\int_\Delta dh_2\, a(h_1)a(h_2)\notag\\
&\qquad\times\big(
e^{-\phi(x_1-x_2)}+
\epsilon u_1(x_1,x_2,h_1,h_2,\epsilon)+\epsilon^2u_2(x_1,x_2,h_1,h_2,\epsilon)\big)\notag\\
&\quad\times
\int_{\Gamma_0}\rho_\mu(d\eta)\big(K^{-1}\big(\epsilon^{-2}
F^{(2)}_{-2}(\cdot,x_1,x_2, h_1, h_2)
+\epsilon^{-1} F^{(2)}_{-1}(\cdot,x_1,x_2, h_1, h_2,\epsilon)\notag\\
&\qquad + F^{(2)}_{0}(\cdot,x_1,x_2, h_1, h_2,\epsilon)\big)  \star
G^{(2)}(\cdot,x_1,x_2, h_1,h_2,\epsilon)\big)(\eta).
\label{ghjsdvfh}
\end{align}

Collecting the coefficients by powers of $\epsilon$, we  get:
\begin{equation}\label{ufyffty} \int_{\Gamma}(H^{(\epsilon )}F)^2(\gamma)\mu(d\gamma)=c_{-2}(\epsilon)\epsilon^{-2}+c_{-1}
(\epsilon)\epsilon^{-1}+c_0(\epsilon)+c_{1}(\epsilon)\epsilon,\end{equation}
where
\begin{align}
c_{-2}(\epsilon)&=\int_{\R^d} z\,dx \int_{\Delta }dh_1\int_\Delta dh_2\,a(h_1) a(h_2)
\int_{\Gamma_0}\rho_\mu(d\eta)
\notag\\
&\quad\times \big(K^{-1}F_{-2}^{(1)}(\cdot,x,h_1,h_2)\star
\big(\big(
e^{-\phi(x-\cdot)}-1
\big)^{\otimes n}\big)_{n=0}^\infty
\big)(\eta)\notag\\
&\quad+\int_{\R^d}z\, dx_1\int_{\R^d}z\,dx_2
\int_\Delta dh_1\int_\Delta dh_2\, e^{-\phi(x_1-x_2)}\notag\\
&\quad\times \int_{\Gamma_0}\rho_\mu(d\eta)
\big(
K^{-1} F_{-2}^{(2)}(\cdot,x_1,x_2,h_1,h_2)
\star\big(\big(
e^{-\phi(x_1-\cdot)-\phi(x_2-\cdot)}-1
\big)^{\otimes n}\big)_{n=0}^\infty\big)(\eta),\notag\\
c_{-1}(\epsilon)&=\int_{\R^d} z\,dx \int_{\Delta }dh_1\int_\Delta dh_2\,a(h_1) a(h_2)
\int_{\Gamma_0}\rho_\mu(d\eta)\Big[
 \big(K^{-1}F_{-2}^{(1)}(\cdot,x,h_1,h_2)\notag\\
&\quad\star \big(n\big(e^{-\phi(x-\cdot)}-1 \big)^{\otimes
(n-1)}\odot g_1^{(1)}(\cdot,x,h_1,h_2,\epsilon)\big)_{n=0}^\infty
\big)(\eta)\notag\\
&\quad+\big(
K^{-1} F_{-1}^{(1)}(\cdot,x,h_1,h_2,\epsilon)\star
\big(\big(
e^{-\phi(x-\cdot)}-1
\big)^{\otimes n}\big)_{n=0}^\infty
\big)(\eta)\Big]\notag\\
&\quad+\int_{\R^d}z\,dx_1 \int_{\R^d} z\,dx_2
\int_\Delta dh_1\int_\Delta dh_2\, a(h_1) a(h_2)\Big[ e^{-\phi(x_1-x_2)}\notag\\
&\quad\times\int_{\Gamma_0}\rho_\mu(d\eta)
\big(K^{-1}F_{-1}^{(2)}(\cdot,x_1,x_2,h_1,h_2,\epsilon)
\star
\big(\big(e^{-\phi(x_1-\cdot)-\phi(x_2-\cdot)}-1
\big)^{\otimes n}
\big)_{n=0}^\infty
\big)(\eta)\notag\\
&\quad + u_1(x_1,x_2,h_1,h_2,\epsilon)
\int_{\Gamma_0}\rho_\mu(d\eta)
\big(K^{-1}F_{-2}^{(2)}(\cdot,x_1,x_2,h_1,h_2)\notag\\
&\quad \star
\big(\big(e^{-\phi(x_1-\cdot)-\phi(x_2-\cdot)}-1
\big)^{\otimes n}
\big)_{n=0}^\infty
\big)(\eta)\notag\\
&\quad+ e^{-\phi(x_1-x_2))}\int_{\Gamma_0}
\big(
K^{-1} F^{(2)}_{-2}(\cdot,x_1,x_2,h_1,h_2)\notag\\
&\quad\star
\big(n\big(e^{-\phi(x_1-\cdot)-\phi(x_2-\cdot)}-1\big)^{\otimes(n-1)}\odot g_1^{(2)}(\cdot,x_1,x_2,h_1,h_2,\epsilon)\big)_{n=0}^\infty
\big)(\eta)\Big],\notag\\
c_0(\epsilon)&=\int_{\R^d} z\,dx\int_\Delta dh_1 \int_\Delta dh_2\,
a(h_1)a(h_2)\int_{\Gamma_0}
\rho_\mu(d\eta)\notag\\
&\quad\times\Big[ \big(
K^{-1} F_{-2}^{(1)}(\cdot,x,h_1,h_2)\notag\\
&\quad\star \big(
n(n-1)(1/2)\big(e^{-\phi(x-\cdot)}-1\big)^{\otimes(n-2)}
\odot \big(g_1^{(1)}(\cdot,x,h_1,h_2,\epsilon)\big)^{\otimes2}\notag\\
&\quad +n\big(e^{-\phi(x-\cdot)}-1\big)^{\otimes(n-1)}\odot
g_2^{(1)} (\cdot,x,h_1,h_2,\epsilon) \big)_{n=0}^\infty
\big)(\eta)\notag\\
&\quad+\big(
K^{-1}F_{-1}^{(1)}(\cdot,x,h_1,h_2,\epsilon)\notag\\
&\quad\star
\big(n\big(e^{-\phi(x-\cdot)}-1\big)^{\otimes(n-1)}\odot g_1^{(1)}
(\cdot,x,h_1,h_2,\epsilon)
\big)_{n=0}^\infty
\big)(\eta)\notag\\
&\quad+\big( K^{-1}F_{0}^{(1)}(\cdot,x,h_1,h_2,\epsilon)\star
\big(\big(e^{-\phi(x-\cdot)}-1\big)^{\otimes n}
\big)_{n=0}^\infty
\big)(\eta)\Big]\notag\\
&\quad+\int_{\R^d}z\,dx_1\int_{\R^d}z\,dx_2
\int_\Delta dh_1\int_\Delta dh_2\, a(h_1)a(h_2)
\Big[ e^{-\phi(x_1-x_2)}
\notag\\
&\quad \times\int_{\Gamma_0}\rho_\mu(d\eta) \big(
K^{-1}F^{(2)}_0(\cdot,x_1,x_2,h_1,h_2,\epsilon)\notag\\
&\quad \star
\big(\big(
e^{-\phi(x_1-\cdot)-\phi(x_2-\cdot)}-1
\big)^{\otimes n}\big)_{n=0}^\infty
\big)(\eta)\notag\\
&\quad+e^{-\phi(x_1-x_2)}\int_{\Gamma_0}\rho_\mu(d\eta) \big(
K^{-1}F^{(2)}_{-1}(\cdot,x_1,x_2,h_1,h_2,\epsilon)\notag\\
&\quad \star
\big(n\big(
e^{-\phi(x_1-\cdot)-\phi(x_2-\cdot)}-1
\big)^{\otimes (n-1)}\odot g_1^{(2)}(\cdot,x_1,x_2,h_1,h_2,\epsilon)\big)_{n=0}^\infty
\big)(\eta)\notag\\
&\quad+e^{-\phi(x_1-x_2)}\int_{\Gamma_0}\rho_\mu(d\eta) \big(
K^{-1}F^{(2)}_{-2}(\cdot,x_1,x_2,h_1,h_2)\notag\\
&\quad \star\big( n(n-1)(1/2)\big(
e^{-\phi(x_1-\cdot)-\phi(x_2-\cdot)}-1
\big)^{\otimes(n-2)}
\odot \big(g_1^{(2)}(\cdot,x_1,x_2,h_1,h_2,\epsilon)\big)^{\otimes 2}
\notag\\
&\quad +n \big(e^{-\phi(x_1-\cdot)-\phi(x_2-\cdot)}-1\big)
^{\otimes(n-1)}\odot g_1^{(2)}(\cdot,x_1,x_2,h_1,h_2,\epsilon)
\big)_{n=0}^\infty \big)
\big)(\eta)\notag\\
&\quad +u_1(x_1,x_2,h_1,h_2,\epsilon) \int_{\Gamma_0}\rho_\mu(d\eta)
\big(
K^{-1}F_{-1}^{(2)}(\cdot,x_1,x_2,h_1,h_2,\epsilon)\notag\\
&\quad \star
\big(
\big(e^{-\phi(x_1-\cdot)-\phi(x_2-\cdot)}-1\big)^{\otimes n}
\big)_{n=0}^\infty
\big)(\eta)\notag\\
&\quad +u_1(x_1,x_2,h_1,h_2,\epsilon) \int_{\Gamma_0}\rho_\mu(d\eta)
\big(
K^{-1}F_{-2}^{(2)}(\cdot,x_1,x_2,h_1,h_2)\notag\\
&\quad \star
\big(n
\big(e^{-\phi(x_1-\cdot)-\phi(x_2-\cdot)}-1\big)
^{\otimes (n-1)}\odot g_1^{(2)}(\cdot,x_1,x_2,h_1,h_2,\epsilon)
\big)_{n=0}^\infty
\big)(\eta)\notag\\
&\quad +u_2(x_1,x_2,h_1,h_2,\epsilon) \int_{\Gamma_0}\rho_\mu(d\eta)
\big(
K^{-1}F_{-2}^{(2)}(\cdot,x_1,x_2,h_1,h_2)\notag\\
&\quad \star \big(
\big(e^{-\phi(x_1-\cdot)-\phi(x_2-\cdot)}-1\big)^{\otimes n}
\big)_{n=0}^\infty \big)(\eta)\Big],\label{yiugytfry}
\end{align}
and $c_1(\epsilon)$ is defined so that  equality
\eqref{ufyffty} holds, i.e., by subtracting from the right hand side of \eqref{ghjsdvfh} the expression
$c_{-2}(\epsilon)\epsilon^{-2}+c_{-1}
(\epsilon)\epsilon^{-1}+c_0(\epsilon)$, given through \eqref{yiugytfry}, and dividing by $\epsilon$.

We evidently have:
\begin{equation}\label{guyf} \int_{\R^d}a(h)h^i\,dh=0,\quad  i\in\{1,\dots,d\},
\end{equation} and therefore
$$c_{-2}(\epsilon)=c_{-1}(\epsilon)=0.$$
Furthermore, as easily  seen $a_1(\epsilon)=O(\epsilon)$ as $\epsilon\to0$.

Below,   we denote $\phi'_i(x):=(\partial/\partial x^i)\phi(x) $ and
$\phi''_{i}(x):=(\partial^2/\partial (x^i)^2)\phi(x)$. So, using
\eqref{guyf}, the equalities
\begin{align*}
&\int_{\R^d} a(h) h^i h^j\,dh=0,\quad i,j\in\{1,\dots,d\},\ i\ne j,
\\ &\int_{\R^d} a(h)(h^i)^2\,dh=c,\quad  i\in\{1,\dots,d\},\end{align*}
  and the dominated convergence theorem, we get
\begin{align}
&\lim_{\epsilon\to0}\int_\Gamma (H^{( \epsilon)}F)^2(\gamma)\,\mu(d\gamma)=\lim_{\epsilon\to0}c_0(\epsilon)\notag\\
&\quad = c^2
\sum_{i,j=1,\dots,d}\bigg[\int_{\R^d}z\,dx\int_{\Gamma_0}\rho_\mu(d\eta)
\Big[ \big(
K^{-1} (\di/\di x^i)F(\cdot\cup x)(\di/\di x^j)F(\cdot\cup x)\notag\\
&\quad\star \big(
n(n-1)\big(e^{-\phi(x-\cdot)}-1\big)^{\otimes(n-2)}\notag\\
&\quad \odot \big( e^{-\phi(x-\cdot)}(-1/2)\phi'_i(x-\cdot)
\big)\odot\big( e^{-\phi(x-\cdot)}(-1/2)\phi'_j(x-\cdot) \big)
\notag\\
&\quad +n\big(e^{-\phi(x-\cdot)}-1\big)^{\otimes(n-1)}\odot
\big( e^{-\phi(x-\cdot)}(-1/2)\phi'_i(x-\cdot)(-1/2)
\phi'_j(x-\cdot) \big) \big)_{n=0}^\infty
\big)(\eta)\notag\\
&\quad+\big( K^{-1} (\di/\di x^i)F(\cdot\cup x)(\di^2/(\di
x^j)^2)F(\cdot\cup x)
\notag\\
&\quad\star
\big(n\big(e^{-\phi(x-\cdot)}-1\big)^{\otimes(n-1)}\odot
\big(e^{-\phi(x-\cdot)}(-1/2)\phi'_i(x-\cdot)\big)
\big)_{n=0}^\infty
\big)(\eta)\notag\\
&\quad+\big( K^{-1} (1/4)(\di^2/(\di x^i)^2)F(\cdot\cup
x)(\di^2/(\di x^j)^2)F(\cdot\cup x) \star
\big(\big(e^{-\phi(x-\cdot)}-1\big)^{\otimes n}
\big)_{n=0}^\infty
\big)(\eta)\Big]\notag\\
&\quad+\int_{\R^d}z\,dx_1\int_{\R^d}z\,dx_2\,\bigg(
 e^{-\phi(x_1-x_2)}
\int_{\Gamma_0}\rho_\mu(d\eta)\notag\\
&\quad\times\Big[ \big( K^{-1} (1/4)(\di^2/(\di x_1^i)^2)F(\cdot\cup
x_1\cup x_2)
(\di^2/(\di x_2^j)^2)F(\cdot\cup x_1\cup x_2)\notag\\
&\quad
\star
\big(\big(
e^{-\phi(x_1-\cdot)-\phi(x_2-\cdot)}-1
\big)^{\otimes n}\big)_{n=0}^\infty
\big)(\eta)\notag\\
&\quad+ \big( K^{-1} (\di/\di x_1^i)F(\cdot\cup x_1\cup x_2)
(\di^2/(\di x_2^j)^2)F(\cdot\cup x_1\cup x_2)
\notag\\
&\quad \star \big(n\big( e^{-\phi(x_1-\cdot)-\phi(x_2-\cdot)}-1
\big)^{\otimes (n-1)}\odot \big(e^{-\phi(x_1-\cdot)
-\phi(x_2-\cdot)}(-1/2)\phi'_i(x_1-\cdot)\big)\big)_{n=0}^\infty
\big)(\eta)\notag\\
&\quad+ \big( K^{-1} (\di/\di x_1^i)F(\cdot\cup x_1\cup x_2)
(\di/\di x_2^j)F(\cdot\cup x_1\cup x_2)
\notag\\
&\quad \star \big( n(n-1)\big(
e^{-\phi(x_1-\cdot)-\phi(x_2-\cdot)}-1
\big)^{\otimes(n-2)}\notag\\
&\quad \odot
\big(e^{-\phi(x_1-\cdot)-\phi(x_2-\cdot)}
(-1/2)\phi_i'(x_1-\cdot)\big)\odot
\big(e^{-\phi(x_1-\cdot)-\phi(x_2-\cdot)}
(-1/2)\phi_j'(x_2-\cdot)\big)
\notag\\
&\quad +n
\big(e^{-\phi(x_1-\cdot)-\phi(x_2-\cdot)}-1\big)^{\otimes(n-1)} \notag\\
&\quad \odot
\big(e^{-\phi(x_1-\cdot)-\phi(x_2-\cdot)}(-1/2)\phi'_i(x_1-\cdot)\phi'_j(x_2-\cdot)\big)
\big)_{n=0}^\infty
\big)(\eta)\Big]\notag\\
&\quad\text{} -(1/2)\phi'_i(x_1-x_2) \int_{\Gamma_0}\rho_\mu(d\eta) \big(
K^{-1}(\di/\di x_1^i)F(\cdot\cup x_1\cup x_2)(\di^2/(\di x_2^j)^2)F(\cdot\cup x_1\cup x_2)\notag\\
&\quad \star
\big(
\big(e^{-\phi(x_1-\cdot)-\phi(x_2-\cdot)}-1\big)^{\otimes n}
\big)_{n=0}^\infty
\big)(\eta)\notag\\
&\quad -\phi'_i(x_1-x_2) \int_{\Gamma_0}\rho_\mu(d\eta) \big(
(\di/\di x_1^i)F(\gamma\cup x_1\cup x_2)(\di/\di x_2^j)F(\gamma\cup
x_1\cup x_2)
\notag\\
&\quad \star \big(n
\big(e^{-\phi(x_1-\cdot)-\phi(x_2-\cdot)}-1\big)^{\otimes
(n-1)}\odot
\big(e^{-\phi(x_1-\cdot)-\phi(x_2-\cdot)}(-1/2)\phi'_j(x_2-\cdot)\big)
\big)_{n=0}^\infty
\big)(\eta)\notag\\
&\quad \text{}-(1/2)\phi'_i(x_1-x_2)(-1/2)\phi'_j(x_2-x_1)\notag\\
&\quad\times \int_{\Gamma_0}\rho_\mu(d\eta) \big( K^{-1} (\di/\di
x_1^i)F(\gamma\cup x_1\cup x_2)
 (\di/\di x_2^j)F(\gamma\cup x_1\cup x_2)
\notag\\
&\quad \star \big(
\big(e^{-\phi(x_1-\cdot)-\phi(x_2-\cdot)}-1\big)^{\otimes n}
\big)_{n=0}^\infty \big)(\eta)\bigg). \label{hgjhftuuy}
\end{align}

Using \eqref{GNZ}, \eqref{jhgftyfwewe}, and \eqref{hjfgytu}, we next have:
\begin{align}
&\int_\Gamma (H^{(\mathrm{dif})}F)^2(\gamma)\mu(d\gamma)=
\frac{c^2}4 \int_\Gamma\mu(d\gamma)\int_{\R^d}z\,dx \exp\big[\la
(-1/2)\phi(x-\cdot),\gamma\ra\big]\notag\\
&\qquad \times\bigg\{
\big(\Delta_xF(\gamma\cup x)\big)^2 -2\Delta_x F(\gamma\cup x)\sum_{u\in\gamma}\la \nabla_x F(\gamma\cup x),\nabla\phi(x-u)\ra\notag\\
&\qquad
+\sum_{u\in\gamma}\la\nabla_x F(\gamma\cup x),
\nabla\phi(x-u)\ra^2
\notag\\&\qquad+
\sum_{u_1\in\gamma}\sum_{u_2\in\gamma\setminus u_1}
\la \nabla_x F(\gamma\cup x),\nabla\phi(x-u_1)\ra
\la \nabla_x F(\gamma\cup x),\nabla\phi(x-u_2)\ra\bigg\}
\notag\\
&\qquad+\frac{c^2}4\int_\Gamma \mu(d\gamma)\int_{\R^d}z\,dx_1\int_{\R^d}z\, dx_2
\exp\big[\la -\phi(x_1-\cdot)-\phi(x_2-\cdot),\gamma\ra\notag\\
&\qquad\text{} -
\phi(x_1-x_2)\big]\bigg\{\Delta_{x_1}F(\gamma\cup x_1\cup x_2)
\Delta_{x_2}F(\gamma\cup x_1\cup x_2)\notag\\
&\qquad- 2\Delta_{x_1}F(\gamma\cup x_1\cup x_2) \sum_{u\in\gamma}
\la\nabla_{x_2}F(\gamma\cup x_1\cup x_2),\nabla\phi(x_2-u)\ra\notag\\
&\qquad -2\Delta_{x_1}F(\gamma\cup x_1\cup x_2) \la
\nabla_{x_2}F(\gamma\cup x_1\cup x_2),\nabla\phi(x_2-x_1)\ra\notag\\
&\qquad+
\sum_{u\in\gamma}\la \nabla_{x_1}F(\gamma\cup x_1\cup x_2),
\phi(x_1-u)\ra
\la \nabla_{x_2}F(\gamma\cup x_1\cup x_2),\nabla
\phi(x_2-u)\ra\notag\\
&\qquad +\sum_{u_1\in\gamma}\sum_{u_2\in\gamma\setminus u_1}
\la  \nabla_{x_1}F(\gamma\cup x_1\cup x_2),\nabla
\phi(x_1-u_1)\ra\,
\la \nabla_{x_2}F(\gamma\cup x_1\cup x_2),
\nabla\phi(x_2-u_2)\ra\notag\\
&\qquad +\sum_{u\in\gamma}2
\la  \nabla_{x_1}F(\gamma\cup x_1\cup x_2),\nabla
\phi(x_1-u)\ra\, 
 \la \nabla_{x_2}F(\gamma\cup x_1\cup x_2),
\nabla\phi(x_2-x_1)\ra\notag\\
&\qquad+ \la  \nabla_{x_1}F(\gamma\cup x_1\cup x_2),\nabla
\phi(x_1-x_2)\ra
\,\la\nabla_{x_2}F(\gamma\cup x_1\cup x_2),\nabla
\phi(x_2-x_1)\ra\bigg\}.\label{ufytfi}
\end{align}

By Lemma \ref{iutgu}, the right hand side of  \eqref{hgjhftuuy} is equal to the right hand side of equality \eqref{ufytfi}. Hence,
\begin{equation}\label{1} \lim_{\epsilon\to 0}\int_\Gamma (H^{(\epsilon)}F)^2(\gamma)\,
\mu(d\gamma) = \int_\Gamma (H^{(
\mathrm{dif})}F)^2(\gamma)\mu(d\gamma).\end{equation} Analogously,
one may also prove that
\begin{equation}\label{2} \lim_{\epsilon\to 0}\int_\Gamma (H^{(\epsilon)}F)(\gamma)
(H^{(\mathrm{dif})}F)(\gamma)\,\mu(d\gamma) = \int_\Gamma (H^{(s,\, \mathrm{dif})}F)^2
(\gamma)\mu(d\gamma).\end{equation}
Now, \eqref{ghdcrtrt} follows from \eqref1 and \eqref2. \quad
$\square$ \vspace{2mm}

\noindent {\it Proof of Corollary }\ref{iuguyfd}. By
Theorem~\ref{uyfytdy}, \cite[Chapter~3, Theorem~3.17]{Dav}, and the
assumption of the corollary,
  we see that, for each $t\ge0$,
$e^{-tH^{( \epsilon)}}\to e^{-tH^{( \mathrm{dif})}}$
strongly in $L^2(\Gamma,\mu)$ as $\epsilon\to0$. To conclude from
here the weak convergence of finite-dimensional distributions, we
proceed  as follows.

We fix any $0\le t_1<t_2<\dots<t_n$, $n\in\mathbb N$. For
$\eps\ge0$, denote by $\mu_{t_1,\dots,t_n}^\eps$ the
finite-dimensional  distribution of the process  $\mathbf
Y^{(s,\,\epsilon)}$ at times $t_1,\dots,t_n$, which is a probability
measure on $\Gamma^n$. Since $\Gamma$ is a Polish space (see e.g.\
\cite{MKM}), by \cite[Chapter II, Theorem~3.2]{Par}, the  measure
$\mu$ is tight on $\Gamma$. Since all the marginal distributions of
the measure $\mu_{t_1,\dots,t_n}^\eps$ are $\mu$, we therefore
conclude that the set $\{\mu_{t_1,\dots,t_n}^\eps\mid\eps>0\}$ is
pre-compact in the space $\mathcal M(\Gamma^n)$ of the probability
measures on $\Gamma^n$ with respect to the weak topology, see e.g.\
\cite[Chapter II, Section~6]{Par}. Hence, the weak convergence of
finite-dimensional distributions follows from the strong convergence
of the semigroups. \quad $\square$ \vspace{2mm}

  \begin{center}
{\bf Acknowledgements}\end{center}

 The authors acknowledge the financial support of the SFB 701 ``Spectral
structures and topological methods in mathematics'', Bielefeld
University and the DFG Project 436 UKR 113/80. We would like to thank
Vladimir Korolyuk, Tobias Kuna, and Michael R\"ockner for useful
discussions.

\end{document}